\documentclass[12pt, a4paper]{article}

\usepackage{a4wide}
\usepackage{amsmath}
\usepackage{amsfonts}
\usepackage{enumerate}
\date{}
\begin{document}
\title{\large\textbf{ON TOPOLOGICAL $\sigma$- IDEALS }}
\author{\normalsize
\textbf{S.BASU \& D.SEN}\\
}
\date{}	
\maketitle
{\small \noindent \textbf{{\textbf{AMS subject classification:}}}} {\small $28A05$, $28A20$, $54C50$ $54E15$.\\}
{\small\noindent \textbf{\textbf{\textmd{\textbf{{Keywords and phrases :}}}}} Topological $\sigma$-ideals, $\mathcal L^{\ast}$-spaces, $\omega$-small system, $\sigma$-algebra admissible with respect to a $\omega$-small system, weakly upper semicontinuous $\omega$-small system, upper semicontinuous small system.}\\
\vspace{.03cm}\\
{\normalsize
\textbf{\textbf{ABSTRACT:}} The concept of $\mathcal S$-topological $\sigma$-ideal in measurable space $(X, \mathcal S)$ was introduced by Hejduk and using a theorem of Wagner on convergence of measurable functions characterized $\mathcal S$-topological $\sigma$-ideals. In this paper, we give a general construction of $\mathcal S$-topological $\sigma$-ideals from structures induced by $\sigma$-algebras and weakly upper semicontinuous $\omega$-small systems. We also show that instead of weak upper semicontinuity, if we use upper semicontinuity, we get $\mathcal S$-uniformizable $\sigma$-ideals. This generalizes the approach of Wagner and Wilczynski metrizing Boolean Lattice of measurable functions.\\

\section{\large{INTRODUCTION}}
\vspace{.1cm}\

\normalsize Much of notations and terminologies used in the sequel are taken from $[12]$ (see also $[11]$)of Wagner and Wilczynski.
\vspace{.01cm}\

Let $\mathcal S$ be a $\sigma$-algebra and $\mathcal I$ be a $\sigma$-ideal of sets in $X$. A property is said to hold $\mathcal I$-almost everywhere (abbr. $\mathcal I$-a.e.) if the set of points in $X$ for which the property does not hold belongs to $\mathcal I$.\\
\vspace{.1cm}

\textbf{DEFINITION $\textbf{1.1}$ :} A sequence $\{f_{n}\}_{n=1}^{\infty}$ of $\mathcal S\Delta\mathcal I$- measurable (i.e. functions measurable with respect to the $\sigma$-algebra generated by $\mathcal S$ and $\mathcal I$) defined on $X$ `converges with respect to $\mathcal I$' to some $\mathcal S\Delta\mathcal I$- measurable function $f$, expressed in writing by $f_{n}\underset{n}{\rightarrow}f$ (w.r.t $\mathcal I$) if every subsequence $\{f_{n_{k}}\}_{
k=1}^{\infty}$ contains a subsequence $\{f_{n_{k_{j}}}\}_{j=1}^{\infty}$ converging to $f$ $\mathcal I$-a.e. This means that the set $\{x\in X: f_{n_{k_{j}}}(x)\not\to f(x)\}\in \mathcal I$.\\
\vspace{.01cm}\

We identify two sets $A$ and $B$ in $\mathcal S\Delta\mathcal I$ if and only if $A\Delta B\in\mathcal I$, thereby obtaining a $\sigma$-complete Boolean algebra $\mathcal S\Delta \mathcal I/\mathcal I$ where the equivalence class containing $A$ is denoted by $[A]$ and for any two classes $[A], [B]\in \mathcal S\Delta\mathcal I/\mathcal I$, the inequality $[A]\leq [B]$ means $A\setminus B\in\mathcal I$. Moreover, $[A]\cup [B]$ means $[A\cup B]$ and ${\displaystyle\bigcup_{n=1}^{\infty}} [A]$ = $\displaystyle {\sup_{n}}[A_{n}]$ means $[{\displaystyle\bigcup_{n=1}^{\infty}}A_{n}]$. In the same manner, we consider two $\mathcal S\Delta\mathcal I$-measurable functions $f$, $g$ on $X$ as equivalent if their difference is equal to zero $\mathcal I$-a.e and in this case denote by $[f]$ the class of all $\mathcal S\Delta\mathcal I$-measurable functions that are equivalent to $f$. The lattice of these classes of equivalences in further denoted by $\mathcal M(\mathcal S\Delta\mathcal I/\mathcal I)$. In this context, it may be said that-\\
\vspace{.01cm}

\textbf{DEFINITION $\textbf{1.2}$ :} A sequence $\{[f_{n}]\}_{n=1}^{\infty}\subseteq \mathcal M(\mathcal S\Delta\mathcal I/\mathcal I)$ `converges with respect to $\mathcal I$' to some $[f]\in \mathcal M(\mathcal S\Delta\mathcal I/\mathcal I)$ expressed in writing by $[f_{n}]\underset{n}{\rightarrow}[f]$ (w.r.t $\mathcal I$) if and only if $f_{n}\underset{n}{\rightarrow}f$ (w.r.t $\mathcal I$).\\
\vspace{.01cm}\

A non empty set $\Omega$ equipped with some notion of convergence is an $\mathcal L^{\ast}$-space if with each sequence $\{p_{n}\}_{n=1}^{\infty}$ convergent in $\Omega$, there is associated a point $p\in\Omega$ such that the following conditions are satisfied:\\
(i) If $\displaystyle\lim_{n \to \infty}p_{n}= p$, then $\displaystyle\lim_{k \to \infty}p_{n_{k}}= p$ for every subsequence $\{p_{n_{k}}\}_{k=1}^{\infty}$ of $\{p_{n}\}_{n=1}^{\infty}$.\\
(ii) If for every $n$, $p_{n}=p$ then $\displaystyle\lim_{n \to \infty}p_{n}= p$.\\
(iii) If a sequence $\{p_{n}\}_{n=1}^{\infty}$ does not converge to $p$, then there exists a subsequence $\{p_{n_{k}}\}_{k=1}^{\infty}$ none of whose subsequences converge to $p$.\\
\vspace{.1cm}\

From Definition $1.1$ and Definition $1.2$, it may be easily observed that the lattice $\mathcal M(\mathcal S\Delta\mathcal I/\mathcal I)$ with the notion of `convergence with respect to the $\sigma$-ideal $\mathcal I$' is an $\mathcal L^{\ast}$-space $[1]$, $[11]$. So it is possible to define an operator $A\rightarrow \overline{A}$ on the set of all subsets of $\mathcal M(\mathcal S\Delta\mathcal I/\mathcal I)$ into itself by setting $\overline{A}=\{[f]: there \hspace{.1cm}exists\hspace{.1cm} a\hspace{.1cm} sequence \hspace{.02cm} \{[f_{n}]\}_{n=1}^{\infty} \subseteq A \hspace{.1cm} such\hspace{.1cm} that\hspace{.05cm} [f_{n}]\underset{n}{\rightarrow}[f]\hspace{.05cm} (w.r.t\hspace{.2cm} \mathcal I)\}$.\\
\vspace{.1cm}\

The operator $A\rightarrow \overline{A}$ satisfies the following properties: $\overline{\emptyset}= \emptyset$, $A\subseteq \overline{A}$, $\overline{A\cup B}= \overline{A}\cup \overline{B}$. In addition, if it also satisfies $\overline{\overline{A}}= \overline{A}$, or, in otherwords,\\
$(*)$ if $[f_{j,n}]\underset{n}{\rightarrow}[f_{j}]$ (w.r.t $\mathcal I$) for each $j$ and $[f_{j}]\underset{j}{\rightarrow}[f]$ (w.r.t $\mathcal I$) implies that there exists subsequences $\{j_{p}\}_{p=1}^{\infty}$ and $\{n_{p}\}_{p=1}^{\infty}$ such that $[f_{j_{p},n_{p}}]\underset{p}{\rightarrow}[f]$ (w.r.t $\mathcal I$),\\
\vspace{.01cm}\

then $\mathcal M(\mathcal S\Delta\mathcal I/\mathcal I)$ becomes a topological space known as Fr$\acute{e}$chet space. There are some examples of $\mathcal L^{\ast}$-spaces which are Fr$\acute{e}$chet spaces and others which are not. It is well-known that the space of Lebesgue measurable functions is a Fr$\acute{e}$chet space whereas the space of functions having Baire property is not so $[11]$. A little uncommon example given by Hejduk $[2]$ shows that if $\mathcal S$ is the $\sigma$-algebra of subsets of a perfect Polish space $X$ containing the family of perfect sets and $\mathcal I$ be any arbitrary totally imperfect $\sigma$-ideal in $X$, then convergence with respect to $\mathcal I$ does not yield a Fr$\acute{e}$chet topology in $\mathcal M(\mathcal S\Delta\mathcal I/\mathcal I)$. Some other examples may be found in $[2]$.\
\vspace{.1cm}\

The following Definition of a `$\mathcal S$-topological ideal' is a slightly modified form of the one already introduced by Hejduk $[1]$.\\
\vspace{.1cm}

\textbf{DEFINITION $\textbf{1.3}$ :} A $\sigma$-ideal $\mathcal I$ on a measurable space $(X,\mathcal S)$ is called a $\mathcal S$-topological ideal if convergence with respect to $\mathcal I$ induces a Fr$\acute{e}$chet topology in $\mathcal M(\mathcal S\Delta\mathcal I/\mathcal I)$.\\
\vspace{.01cm}\

In $[12]$, Wagner and Wilczynski studied the problem of topologizing a space of measurable functions. They proved that-\\
\vspace{.1cm}

\textbf{THEOREM $\textbf{1.4}$ :} If $\mathcal S$ is a $\sigma$-algebra and $\mathcal I (\subseteq \mathcal S)$ a $\sigma$-ideal of sets in $X$, then convergence with respect to $\mathcal I$ yields a Fr$\acute{e}$chet topology in $\mathcal M(\mathcal S/\mathcal I)$ if and only if the $\sigma$-complete Boolean algebra $\mathcal S/\mathcal I$ satisfies the following conditions:\\
(*) for every $[D]\in\mathcal S/\mathcal I$, $[D]\neq\mathcal I$ and for every double sequence $\{[B_{j,n}]\}_{j,n=1}^{\infty} (\subseteq \mathcal S/\mathcal I)$ such that-\\
(i) $[B_{j,n}]\leq [B_{j,n+1}]$ for every $j,n\in \mathbb{N}$\\
(ii) $[\displaystyle\bigcup _{n=1}^{\infty}B_{j,n}]=[D]$ for every $j\in \mathbb{N}$\\
\vspace{.1cm}\

there exists a positive integer valued functions $n(i,j)$ such that $[\displaystyle\bigcup_{n=1}^{\infty}\bigcap_{n=1}^{\infty} B_{j,n}(i,j)]=[D]$.\
\vspace{.1cm}\

Based on the above result of Wagner and Wilczynski, we may rewrite the characterization of $\mathcal S$-topological ideal by Hejduk in the present settings-\\
\vspace{.1cm}

\textbf{THEOREM $\textbf{1.5}$ :} A $\sigma$-ideal $\mathcal I$ in a measurable space $(X,\mathcal S)$ is a $\mathcal S$-topological ideal ($\mathcal I$ is not necessarily contained in $\mathcal S$) if and only if the $\sigma$-complete Boolean algebra $\mathcal S\Delta \mathcal I/\mathcal I$ satisfies condition (*).\\
\vspace{.1cm}\

In this paper, we introduce a general process of constructing $\mathcal S$-topological $\sigma$-ideal on a non-empty set $X$ using $\sigma$-algebra $\mathcal S$ and weakly upper semi continuous small system which is admissible with respect to the $\sigma$-algebra $\mathcal S$. If weak upper semicontinuity is replaced by upper semicontinuity, the $\sigma$-ideal turns out to be uniformizable. Such types of structures were earlier used by authors like Niewiarowski $[6]$, Hejduk and Wajch $[4]$ in connection with generalization of Fr$\acute{e}$chet's theorem characterizing compactness in the family of measurable functions in the sense of convergence with respect to finite measure. \\
\vspace{.1cm}

\section{\large{PRELIMINARIES AND RESULTS}}

\vspace{.1cm}\
\normalsize\hspace{.2cm} Let $X$ be a non empty set. By\\
\vspace{.1cm}

\textbf{DEFINITION $\textbf{2.1}$ :} A $\omega$-small system on $X$, we mean a sequence $\{\mathcal N_{n}\}_{n=1}^{\infty}$ where $\mathcal N_{n}$ is class of subsets of $X$ satisfying the following set of conditions:\\
(i) $\emptyset\in \mathcal N_{n}$ for all $n$.\\
(ii) If $E\in\mathcal N_{n}$ and $F\subset E$, then $F\in\mathcal N_{n}$; i.e, each $\mathcal N_{n}$ is a hereditary class.\\
(iii) $E\in \mathcal N_{n}$ and $F\in{\displaystyle{\bigcap_{n=1}^{\infty}}\hspace{.01cm}{\mathcal N_{n}}}$ implies that $E\cup F\in\mathcal N_{n}$.\\
(iv) For each $m$, there exists $m^{\prime}> m$ such that for any one-to-one correspondence $k\rightarrow n_{_{k}}$ with $n_{k} > m^{'}$, $\displaystyle{\bigcup_{k=1}^{\infty}}{E_{n_{k}}}\in \mathcal N_{_{m}}$ whenever $E_{n{_{k}}}\in \mathcal N_{n{_{k}}}$.\\
(v) For any $m$, $n$, there exists $p> m,n$ such that $\mathcal N_{p}\subseteq \mathcal N_{m}$ and $\mathcal N_{p}\subseteq\mathcal N_{n}$; i.e, $\{\mathcal N_{n}\}_{n=1}^{\infty}$ is directed.\\
\vspace{.1cm}\

We set $\mathcal N_{\infty}=\displaystyle\bigcap_{n=1}^{\infty}{\mathcal N_{n}}$. From conditions (i), (ii) and (iv) of Definitions $2.1$, it follows that $\mathcal N_{\infty}$ is $\sigma$-ideal.\\
\vspace{.1cm}\

The above Definition is a modified version of the notion of `small system' originally introduced by Ri$\acute{c}$an $[7]$, $[8]$, Ri$\acute{c}$an and Neubrunn $[9]$ and used by others $[5]$, $[10]$ including those whose names have been mentioned at the begining of the introduction.\\
\vspace{.1cm}

\textbf{DEFINITION $\textbf{2.2}$ :} A $\sigma$-algebra $\mathcal S$ on $X$ is admissible with respect to a $\omega$-small system $\{\mathcal N_{\infty}\}_{n=1}^{\infty}$ on $X$ if for each $n$,\\
(i) $\mathcal S\setminus\mathcal N_{n}\neq\emptyset\neq \mathcal S\cap\mathcal N_{n}$.\\
(ii) Each $\mathcal N_{n}$ has a $\mathcal S$-base; i.e, each $E\in\mathcal N_{n}$ is contained in some $F\in \mathcal S\cap\hspace{.02cm}\mathcal N_{n}$. In otherwords, we may express this by saying that $\{\mathcal N_{\infty}\}_{n=1}^{\infty}$ is $\mathcal S$-regular.\\
(iii) $\mathcal S\setminus \mathcal N_{n}$ satiesfies countable chain condition (c.c.c) which means that any arbitrary collection of mutually disjoint sets from $\mathcal S\setminus \mathcal N_{n}$ is atmost countable.\\
\vspace{.01cm}\
we further add that- \\
\vspace{.1cm}

\textbf{DEFINITION $\textbf{2.3}$ :} A $\omega$-small system $\{\mathcal N_{n}\}_{n=1}^{\infty}$ is weakly upper semicontinuous relative to a $\sigma$-algebra $\mathcal S$ on $X$ if for every nested sequence $\{E_{n}\}_{n=1}^{\infty}$ $(E_{n+1}\subseteq E_{n})$ of sets from $\mathcal S$ satisfying $\displaystyle\bigcap_{n=1}^{\infty}{E_{n}}=\emptyset$, there exists for every $k$ some $n_{k}$ such that $E_{n_{k}}\in\mathcal N_{k}$;\\
\vspace{.1cm}\
and\\
\vspace{.1cm}

\textbf{DEFINITION $\textbf{2.4}[6]$ :} A $\omega$-small system $\{\mathcal N_{n}\}_{n=1}^{\infty}$ is  upper semicontinuous relative to a $\sigma$-algebra $\mathcal S$ on $X$ if for every nested sequence $\{E_{n}\}_{n=1}^{\infty}$ $(E_{n+1}\subseteq E_{n})$ of sets from $\mathcal S$, $E_{n}\notin \mathcal N_{m}$ for some $m$ and $n=1,2,\ldots$ implies that $\displaystyle\bigcap_{n=1}^{\infty}{E_{n}}\notin \mathcal N_{\infty}$.\\
\vspace{.01cm}\

By virtue of the condition (i) of Definition $2.1$, one can readily see that upper semi continuity is a stronger notion than weak upper semi continuity.\\
\vspace{.01cm}\

In none of the above Definitions, it is assumed that $\mathcal N_{n}$ is a subclass of $\mathcal S$ for any $n$. So in general $\mathcal N_{\infty}$ is not a subclass of $\mathcal S$. However, we can consider the $\sigma$-algebra generated by $\mathcal S$ and $\mathcal N_{\infty}$ and denote this by $\widetilde{\mathcal S}$. In fact, every element of $\widetilde{\mathcal S}$ admits a representation of the form $(A\setminus B)\cup C$ where $A\in\mathcal S$ and $B,C (\in\mathcal N_{\infty})$. If $\mathcal S$ is admissible with respect to $\{\mathcal N_{n}\}_{n=1}^{\infty}$, then from conditions (ii) and (iii) of Definition $2.2$, it follows that $\widetilde{\mathcal S}\setminus \mathcal N_{\infty}$ satisfies countable chain condition. This in otherwords means that $\widetilde{\mathcal S}$ satisfies countable chain condition (c.c.c).\\
\vspace{.01cm}\

Now to show that $\mathcal N_{\infty}$ is a $\mathcal S$-topological ideal we proceed in steps. We first prove the following theorem\\
\vspace{.01cm}

\textbf{THEOREM $\textbf{2.5}$ :} Let $\mathcal S$ be a $\sigma$-algebra admissible with respect to a $\omega$-small system $\{\mathcal N_{n}\}_{n=1}^{\infty}$ which is weakly upper semicontinuous relative to $\mathcal S$. Let $E\in \mathcal S\setminus \mathcal N_{\infty}$ and $\{E_{j,n}\}_{j,n=1}^{\infty}$ be a double sequence of sets from $\mathcal S$ satisfying:\\
(i) $E_{j,n}\subseteq E_{j,n+1}$ for every $j,n =1,2,\ldots$\\
(ii) $\displaystyle\bigcup_{n=1}^{\infty} {E_{j,n}}= E$ for every $j=1,2,\ldots$\\
\vspace{.01cm}\

Then there exists a sequence $\{n_{j}\}_{j=1}^{\infty}$ of naturals numbers such that $\displaystyle\bigcap_{j=1}^{\infty}{E_{j,n_{j}}}\notin \mathcal N_{\infty}$.\\
\vspace{.1cm}

\textbf{PROOF :} We observe that there exists a natural number $m$ such that no subset $M$ of $E$ and its complement $E\setminus M$ (in $E$) can both belong to $\mathcal N_{m}$. Since $E\notin \mathcal N_{\infty}$, there exists a natural number $r$ such that $E\in \mathcal S\setminus \mathcal N_{r}$ and also because there exist $p,q>r$ (by (iv), Definition $2.1$) such that $\mathcal N_{p}\cup\mathcal N_{q}\subseteq \mathcal N_{r}$ and $m>p,q$ (by (v), Definition $2.1$) such that $\mathcal N_{m}\subseteq \mathcal N_{p}$, $\mathcal N_{m}\subseteq \mathcal N_{q}$, it will follow that if both $M$ and $E\setminus M$ belong to $\mathcal N_{m}$ then their union which is $E$ should be in $\mathcal N_{r}$ - a contradiction.\\
\vspace{.01cm}\

We set $G_{j,n}=E\setminus E_{j,n}$ , $j=1,2,\ldots$, then $\displaystyle\bigcap_{n=1}^{\infty} {E_{j,n}}= \emptyset$ for each $j=1,2,\ldots$. Since $\{\mathcal N_{n}\}_{n=1}^{\infty}$ is weakly upper semicontinuous relative to $\mathcal S$, for each $j$, there exists $n_{j}$ such that $G_{j,n_{j}}\in \mathcal N_{n_{j}}$ and by virtue of condition (iv) of Definition $2.1$, we may choose $n_{j}>m^{'}$ so that $\displaystyle\bigcup_{j=1}^{\infty}{G_{j,n_{j}}}\in \mathcal N_{m}$. It then follows from the above observation that $\displaystyle\bigcap_{j=1}^{\infty}{E_{j,n_{j}}}\notin \mathcal N_{m}$. Consequently, $\displaystyle\bigcap_{j=1}^{\infty}{E_{j,n_{j}}}\notin \mathcal N_{\infty}$.\\
\vspace{.1cm}\

Using a result of Hejduk (Proposition $1$, $[3]$) we may now write that - \\
\vspace{.1cm}

\textbf{THEOREM $\textbf{2.6}$ :} If a $\sigma$-algebra $\mathcal S$ on $X$ is admissible with respect to a $\omega$-small system $\{\mathcal N_{\infty}\}_{n=1}^{\infty}$ which is also weakly upper semicontinuous relative to $\mathcal S$; $[E]\in\widetilde{\mathcal S}/\mathcal N_{\infty}$ and $\{[E_{j,n}]\}_{j,n=1}^{\infty}$ ($\subseteq \widetilde{\mathcal S}/\mathcal N_{\infty}$) is a double sequence satisfying:\\
(i) $[E_{j,n}]\leq [E_{j,n+1}]$ for each $j,n =1,2,\ldots$\\
(ii) $[\displaystyle\bigcup_{n=1}^{\infty} {E_{j,n}}]= [E]$ for every $j=1,2,\ldots$\\
\vspace{.01cm}\

Then there exists a sequence $\{n_{j}\}_{j=1}^{\infty}$ of naturals numbers such that $\displaystyle\bigcap_{j=1}^{\infty}{[E_{j,n_{j}}]}\notin \mathcal N_{\infty}$.\\
\vspace{.01cm}\
Since $\widetilde{\mathcal S}\setminus \mathcal N_{\infty}$ satisfies countable chain condition, there exists a positive integer valued function $n(i,j)$ of positive integers $i$, $j$ such that $[\displaystyle \bigcup_{i=1}^{\infty}\bigcap_{j=1}^{\infty}{E_{j,n(i,j)}}]= [E]$\\
\vspace{.01cm}\
\hspace{.3cm} Hence under the hypothesis of the above theorem, condition (*) (Theorem $1.4$) is fulfilled for the Boolean lattice $\widetilde{\mathcal S}/\mathcal N_{\infty}$. So from the Hejduk's theorem (Theorem $1.4$), it follows that\\
\vspace{.1cm}

\textbf{THEOREM $\textbf{2.7}$ :} The $\sigma$-ideal $\mathcal N_{\infty}$ is a $\mathcal S$-topological ideal on $X$.\\
\vspace{.1cm}\

We next show that the $\sigma$-ideal $\mathcal N_{\infty}$ is uniformizable in the sense that we can define an uniformity on $\mathcal M(\widetilde {S}/\mathcal N_{\infty})$ if instead of weak upper semicontinuity, we use the stronger version of upper semicontinuity. This generalizes the approach of Wagner and Wilczynski $[12]$ metrizing  Boolean lattice of measurable functions in abstract settings.\\
\vspace{.1cm}\

For any $r>0$ and natural numbers $n$, we set\\ $V_{r,n}=\{([f], [g])\in\mathcal M(\widetilde{\mathcal S}/\mathcal N_{\infty})\times\mathcal M(\widetilde{\mathcal S}/\mathcal N_{\infty}):\{x\in X:\mid{f(x)-g(x)}\mid\geq r\}\in\mathcal N_{n}\}$. We claim that the family $\{V_{r,n}: r>0, n=1,2,\ldots\}$ constitutes a base for some uniformity in $\mathcal M(\widetilde{\mathcal S}/\mathcal N_{\infty})$. Evidently, by condition (i) of Definition $2.1$, the diagonal $\Delta\subseteq V_{r,n}$ for every $r>0$ and $n$. Also $V_{r,n}^{-1}= V_{r,n}$ for every $r>0$ and $n$. Now for any natural number $n$, it is possible to choose by condition (iv) of Definition $2.1$, natural numbers $m$ and $p$ such that $\mathcal N_{m}\cup\mathcal N_{p}\subseteq \mathcal N_{n}$. Also, as $\{x\in X: \mid f(x)-h(x)\mid\geq r\}\subseteq \{x\in X: \mid f(x)-g(x)\mid\geq r/2\}\cup \{x\in X: \mid g(x)-h(x)\mid\geq r/2\}$ for every $r>0$, therefore $V_{r/2,p}\circ V_{r/2,m}\subseteq V_{r,n}$. Finally, since for any two natural numbers $m,n$ there exists a natural number $p$ such that $\mathcal N_{p}\subseteq\mathcal N_{m}, \mathcal N_{n}$ by condition (v) of Definition $2.1$, therefore $V_{r,p}\subseteq V_{r,m}\cap V_{r,n}$.\\
\vspace{.1cm}\

We now show that `` convergence with respect to the uniformity induced by the family $\{V_{r,n}: r>0, n=1,2,\ldots\}$ " is equivalent to the `` convergence with respect to the $\sigma$-ideal $\mathcal N_{\infty}$ "\\
\vspace{.1cm}\

Let $\{[f_{k}]\}_{k=1}^{\infty}\subseteq \mathcal M(\widetilde{\mathcal S}/\mathcal N_{\infty})$ converges to $[f]$ in the above uniformity. Then for every subsequence $\{f_{k_{j}}\}_{j=1}^{\infty}$, there exists a subsequence $\{f_{k_{j_{_{n}}}}\}_{n=1}^{\infty}$ such that $([f_{k_{j_{_{n}}}}], [f])\in V_{\frac{1}{2^{k_{j_{_{n}}}}},k_{j_{_{n}}}}$, where the subsequence can be so chosen that $\displaystyle\bigcap_{m=1}^{\infty}\bigcup_{n=m}^{\infty}{\mathcal N_{k_{j_{_{n}}}}}\subseteq\displaystyle\bigcap_{m=1}^{\infty}{\mathcal N_{m}}= \mathcal N_{\infty}$. This choice is possible by repeated application of condition (iv) of Definition $2.1$. Hence $[f_{n}]\underset{n}{\rightarrow}[f]$ (w.r.t $\mathcal N_{\infty}$).\\
\vspace{.1cm}\

Conversely, suppose the $\{[f_{n}]\}_{n=1}^{\infty}$ does not converge to $[f]$ in the above uniformity. This means that there exist some $r_{0}>0$, a natural number $n_{0}$ and a subsequence $\{f_{n_{_{k}}}\}_{k=1}^{\infty}$ such that $\{x\in X: \mid f_{n_{_{k}}}(x)-f(x)\mid\geq r_{0}\}\notin \mathcal N_{n_{_{0}}}$.\\ We put $E_{k}=\{x\in X : x\in X: \mid f_{n_{_{k}}}(x)-f(x)\mid\geq r_{0}\}$ and so $F_{m}=\displaystyle\bigcup_{k\geq m}{E_{k}} \notin \mathcal N_{n_{_{0}}}$ for every $m$ by condition (ii) of Definition $2.1$. But $F_{m}= G_{m}\Delta H_{m}$, where $G_{m}\in\mathcal S$ and $\mathcal H_{m}\in\mathcal N_{\infty}$. Since $\mathcal S$ is admissible with respect to $\{\mathcal N_{n}\}_{n=1}^{\infty}$, by (ii) of Definition $2.2$, $H_{m}\subseteq K_{m}\in \mathcal N_{\infty}\cap \mathcal S$ for every $m$. Therefore, $F_{m}\supseteq G_{m}\setminus\displaystyle\bigcup_{j=1}^{m}{K_{j}}\notin \mathcal N_{n_{_{0}}}$ for otherwise, by (iii) of Definition $2.1$, $F_{m}\in \mathcal N
_{\infty}$. But $G_{m}\setminus\displaystyle\bigcup_{j=1}^{m}{K_{j}}\in\mathcal S$ and so by upper semicontinuity of $\{\mathcal N_{n}\}_{n=1}^{\infty}$ relative to $\mathcal S$, $\displaystyle\bigcap_{n=1}^{\infty}(G_{n}\setminus\displaystyle\bigcup_{j=1}^{n}{K_{j}})\notin \mathcal N_{\infty}$. Consequently, by (ii) of Definition $2.1$, $\displaystyle\bigcap_{n=1}^{\infty}{F_{n}}\notin \mathcal N_{\infty}$. Hence $[f_{n}]\not\to _{n}[f]$ (w.r.t $\mathcal N_{\infty}$).\\
\vspace{.1cm}\

Thus from what we have deduced above, it follows that\\
\vspace{.1cm}

\textbf{THEOREM $\textbf{2.8}$ :} If $\mathcal S$ is a $\sigma$-algebra on $X$ admissible with respect to a $\omega$-system $\{\mathcal N_{n}\}_{n=1}^{\infty}$ which is upper semicontinuous relative to $\mathcal S$, then the $\sigma$-ideal $\mathcal N_{\infty}$ is uniformizable.\\
\vspace{.1cm}

{\textbf{REMARKS:}} We observed in $[12]$, that condition (*) does not imply countable chain condition (c.c.c). But in our case, in order to show that $\widetilde{\mathcal S}/\mathcal N_{\infty}$ satisfies (*) (Theorem $2.6$), we need countable chain condition in $\widetilde{\mathcal S}\setminus\mathcal N_{\infty}$ as a part of our hypothesis. However, countable chain condition may be relaxed by virtue of upper semicontinuity. In fact, in proving Theorem $2.8$, we do not even need the full force of admissibility.\\
\vspace{.08cm}
\hspace{1cm}

\vspace{.1cm}

\underline{\textbf{AUTHOR'S ADDRESS}}\\

{\normalsize
\textbf{S.Basu}\\
\vspace{.02cm}
\hspace{.35cm}
\textbf{Dept of Mathematics}\\
\vspace{.02cm}
\hspace{.35cm}
\textbf{Bethune College , Kolkata} \\
\vspace{.01cm}
\hspace{.35cm}
\textbf{W.B. India}\\
\vspace{.02cm}
\hspace{.1cm}
\hspace{.29cm}\textbf{{e-mail : sanjibbasu08@gmail.com}}\\
%\vspace{.0001cm}

\textbf{D.Sen}\\
\vspace{.02cm}
\hspace{.35cm}
\textbf{Saptagram Adarsha vidyapith (High), Habra , $24$ Parganas (N)} \\
\vspace{.01cm}
\hspace{.35cm}
\textbf{W.B. India}\\
\vspace{.02cm}
\hspace{.1cm}
\hspace{.29cm}\textbf{{e-mail : reachtodebasish@gmail.com}}
}

\end{document}